\documentclass[11pt]{article}
\newcommand{\proof}{\noindent {\bf Proof: }}

\newtheorem{theorem}{Theorem}

\newtheorem{proposition}{Proposition}
\newtheorem{defi}{Definition}
\def\qed{\hfill $\Box$}

\usepackage{amssymb}
\usepackage{graphicx}

\begin{document}

\title{Conics in normed planes}
\author{\'{A}kos G. Horv\'{a}th\footnote{Department of Geometry, Budapest University of Technology and Economics
1521 Budapest, Hungary}\,\,\, and \, Horst Martini\footnote{Faculty of Mathematics, Chemnitz University of Technology, 09107 Chemnitz, Germany}}
\date{~}

\maketitle

\begin{abstract}
We study the generalized analogues of conics for normed planes by using the
following natural approach: It is well known that there are different
metrical definitions of
conics in the Euclidean plane.
We investigate how these definitions extend to normed planes, and we
show that in this more general framework these different definitions
yield, in almost all cases, different classes of curves.
\end{abstract}



\section{Introduction}

 We present a systematic investigation of possible definitions of conics extended to normed (or Minkowski) planes. In the Euclidean situation the metric definitions of conics and the analytic one,
 namely defining them as family of curves of second order, clearly yield the same type of curves; so we have various different definitions
 of the same class of curves. In normed planes neither the metric definitions nor the analytic one yield the same type of curves. Furthermore, it is not clear what the notions ``curve of second order'', ``cone of second order''
 or ``sections of a cone'' mean. We consider the usual metric definitions of conics in the Euclidean plane, adopt them for normed planes and list
 various properties of the resulting classes of curves.\\
  By $X$ we denote a \emph{normed} or \emph{Minkowski plane}, i.e., the affine
plane equipped with a norm $\| \cdot \|$ determined by the \emph{unit ball}
$$
K = \{ {\bf x}\in X: \|{\bf x}\| \le 1\},
$$
which is a compact, convex set centered at the \emph{origin} ${\bf o}$, being
from the interior of $K$. The boundary of $K$ is the \emph{unit circle} $S$
of $X$.
We say that $K$ (or the norm induced by $K$) is \emph{strictly convex} if
$S$ contains no proper segment. We
use small letters like ${\bf x}, {\bf y}$ for points/vectors in $X$, and the
symbol $[{\bf x} , {\bf y}]$ describes the closed segment with different endpoints
${\bf x}$ and ${\bf y}$.
It is clear that the notion of \emph{bisector}
$Bis({\bf x} , {\bf y})$ of two different points $x,y$ from
$X$, defined by
$$
Bis({\bf x}, {\bf y}) = \{ {\bf z} \in X: \| {\bf x} - {\bf z}\| = \| {\bf z} - {\bf y}\|\}
$$
can be suitably extended to bisectors of two point sets
(instead of the two points ${\bf x}$ and ${\bf y}$); see [4] and [5] for that
notion. The \emph{d-segment} with different endpoints ${\bf x}$ and ${\bf y}$ from $X$
is the point set defined by
$$
[{\bf x}, {\bf y}]_d = \{ {\bf z}\in X :\| {\bf x} -{\bf z}\| +\| {\bf z}- {\bf y}\| = \| {\bf x}- {\bf y}\| \};
$$
see $\S 9$ in [3]. Of course, depending on the shape of $K$, bisectors
and $d$-segments are geometrically interesting and can even be
two-dimensional; see again [4], [5], and $\S 9$ in [3]. But it
should be noticed that, for strictly convex norms, any bisector is
homeomorphic to a line and any $d$-segment is a usual segment.

\section{Ellipses defined by metric properties}

First we consider the usual metric definitions of ellipses in the Euclidean plane and examine their analogues (i.e., generalizations) for normed planes.

\begin{figure}[ht]
\centering
  \includegraphics[scale=0.5]{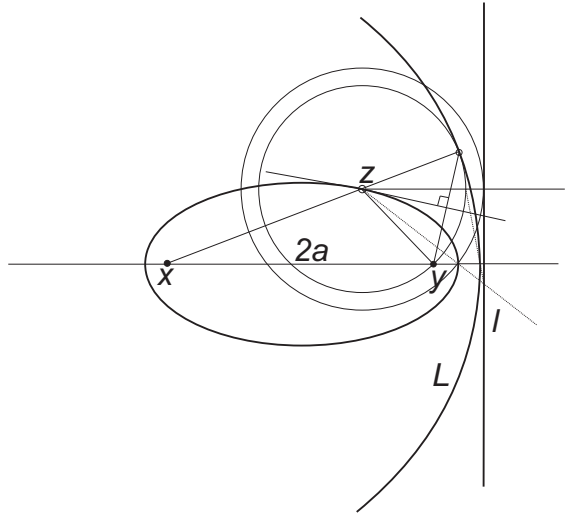}\\
  \caption{ An ellipse in the Euclidean case (Definition 3)}
\end{figure}

For this purpose we give the basic figure of an ellipse in the Euclidean plane (see Fig.1), containing the foci ${\bf x}$ and ${\bf y}$, a (variable) point ${\bf z}$ corresponding to the distance sum $|{\bf z}-{\bf x}|+|{\bf z}-{\bf y}|=2a$ (with $| \cdot |$ denoting the Euclidean norm), the tangent line at this point ${\bf z}$, the reflected image of ${\bf y}$ in this line, the leading circle $L$ which is the locus of such reflected images, and the leading line $l$, defined as the common tangent of the circles having radius $\frac{a}{c}|{\bf z}-{\bf y}|$, where $c$ is half the distance of the foci.

In normed planes we have three different possibilities to define ellipses metrically. Until now, only the first one was investigated (see \cite{wu}).
So the following three definitions refer to a normed plane $X$.

\begin{defi}[based on foci]
Let ${\bf x},{\bf y}\in X$, ${\bf x}\not ={\bf y}$, and $2a\geq 2c=\|{\bf x}-{\bf y}\|$.  The set
$$
E({\bf x},{\bf y}, a)=\{{\bf z} \in X: ||{\bf z}-{\bf x}|| +||{\bf z}-{\bf y}||=2a \}
$$
is called \emph{the ellipse defined by its foci} $x$ and $y$.
\end{defi}

\begin{defi}[based on a leading circle and one focus]
Let $L:=(2a)\cdot K$ be a homothetic copy of the unit disk $K$, and ${\bf x}\in L$ be an arbitrary point from it. The locus of points ${\bf z}\in X$ for which there is a positive $\varepsilon$ such that ${\bf z}+\varepsilon K$ touches $L$ and contains ${\bf x }$ on its boundary is called \emph{the ellipse defined by its leading circle and its focus} {\bf x}.
\end{defi}

\begin{defi}[based on a leading line and a focus]
Let $l$ be a straight line, ${\bf x}$ a point, and $\gamma =\frac{a}{c}$ a ratio larger than 1. The locus of points ${\bf z}\in X$, for which there is a positive ${\varepsilon}$ such that the boundary of the disk ${\bf z}+\varepsilon K$ contains {\bf x} and the disk ${\bf z}+\gamma (\varepsilon K)$ touches the line $l$, is called \emph{the ellipse defined by its leading line and its focus} {\bf x}.
\end{defi}

    The equivalence of these definitions for the Euclidean subcase is well known, and this can be easily checked  with the help of Fig. 1. We will prove that, while the first two definitions are  equivalent also in normed planes, the third one yields a basically different class of curves.

\begin{proposition}
 In any normed plane the following holds: an ellipse, defined by its foci, is always an ellipse defined by its leading circle and a focus, and the converse statement is also true. On the other hand, an ellipse defined by its leading line and a focus is not necessarily an ellipse defined by its foci, and again the converse is true.
\end{proposition}

\proof  First we consider an ellipse which is defined by its leading circle. Thus we have two disks $D_{\bf z}={\bf z}+\varepsilon K$ and $L=(2a)\cdot K$ in touching position. Then the line joining their centers contains a point from the intersection of their boundaries. Call this point ${\bf p}$. Then
$$
\varepsilon=\|{\bf p}-{\bf z}\|=\|{\bf z}-{\bf x}\|\,.
$$
Thus
$$
2a=\|{\bf z}\|+\|{\bf z}-{\bf x}\|,
$$
implying that ${\bf z}\in E({\bf 0},{\bf x}, a)$.

\begin{figure}[ht]
  \centering
  \includegraphics[scale=0.5]{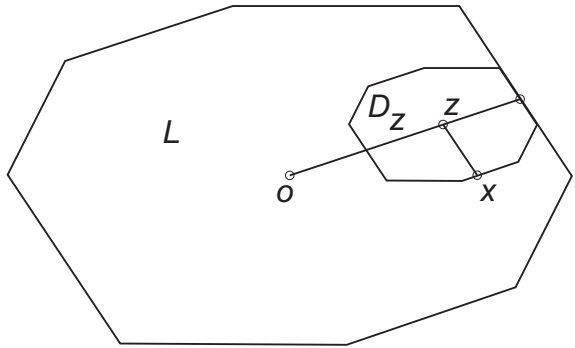}\\
  \caption{ The equivalence of Definitions 1 and 2}
\end{figure}

On the other hand, with the same notation for a point ${\bf z}\in E({\bf 0},{\bf x}, a)$,  its definition yields $\|{\bf z}\|<2a$. Consider the disk $D_{\bf z}={\bf z}+(2a-\|{\bf z}\|)K$. This disk is touching $L$. Since
$$
2a=\|{\bf z}\|+\|{\bf z}-{\bf x}\|
\mbox{ and so }
\|{\bf z}-{\bf x}\|=(2a-\|{\bf z}\|),
$$
${\bf x}$ is on the boundary of $D_{\bf z}$.

 Now we give examples showing that there exists an ellipse which is defined via its leading line but is not an ellipse defined via its foci, and conversely.

In Fig. 3 we can see
that there is an ellipse following the third definition which is not centrally symmetric. By Theorem 2 of \cite{wu} it is not an ellipse by the first definition. In our example the norm is the $L_\infty$ norm, and the leading line $l$ and the focus ${\bf x}$ are in ``symmetric position'' with respect to the circle of this Minkowski plane, which is a square.

\begin{figure}[ht]
 \centering
  \includegraphics[scale=0.5]{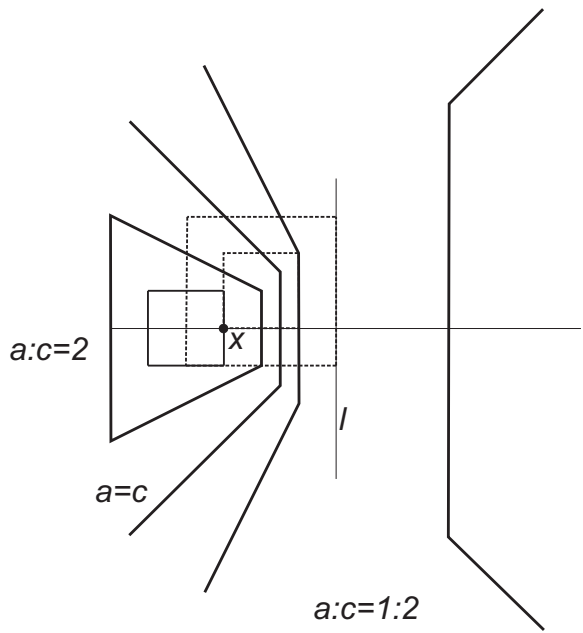}\\
  \caption{Conics on the $l_\infty$ plane}
\end{figure}

\begin{figure}[ht]
    \centering
  \includegraphics[scale=0.5]{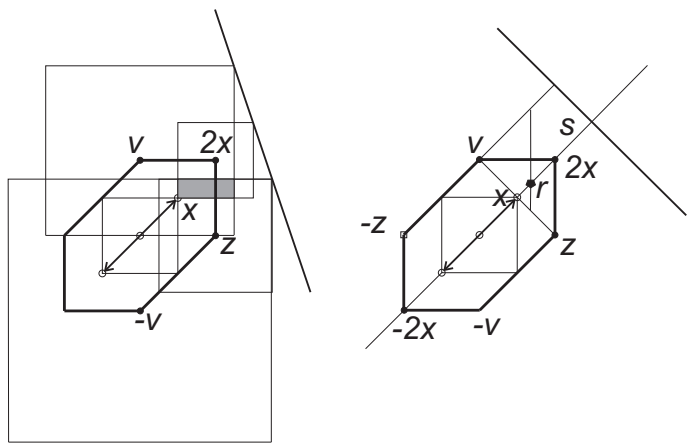}\\
  \caption{A metric ellipse which has no leading line}
\end{figure}

Conversely, consider the ellipse $E(-{\bf x},{\bf x}, 2)$ defined by its foci and shown in Fig. 4.
First we can see that if it is also an ellipse defined by its leading line, then the leading line $l$ and the new focus ${\bf x}'$ have to be in ``symmetric position'' with respect to the line joining the original foci. ``Symmetric'' means that this line is parallel to a diagonal of the unit square. In fact, if this is not the case, we get a figure as shown on the left side of Fig. 4. The squares $S_{2{\bf x}}, S_{\bf v}, S_{\bf z}, S_{-{\bf v}}$ with centers $2{\bf x}$, ${\bf v}$, ${\bf z}$, $-{\bf v}$, respectively, touch $l$. The focus has to lie in the shaded rectangle, as the common point of the boundaries of homothetic copies $2{\bf x}+\frac{c}{a}S_{2{\bf x}}$, ${\bf v}+\frac{c}{a}S_{{\bf v}}$ and ${\bf z}+\frac{c}{a}S_{{\bf z}}$ of such squares (with a homothety ratio smaller than 1). On the other hand, the boundary of the square $-{\bf v}+\frac{c}{a}S_{-{\bf v}}$ intersects the shaded rectangle in a segment parallel to that one in which it is intersected by ${\bf z}+\frac{c}{a}S_{{\bf z}}$. So it is

We now assume that $l$ and ${\bf x}'$ have symmetric position (see the right side of Fig. 4).
If this holds and the Euclidean distance of $l$ and $2{\bf x}$ is $s$, and that of ${\bf x'}$ and ${\bf x}$ is $r$, then, using the fact that the points $2{\bf x}$, $-2{\bf x}$ and ${\bf v}$
have to lie on the new ellipse, we have the equalities
$$
\frac{r}{s}=\frac{4-r}{4+s}=\frac{2-r}{1+s}\,,
$$
implying that
$$
s=1 \mbox{ and } r=\frac{2}{3}
$$
and showing that $\frac{a}{c}=\frac{2}{3}$.
Thus the leading line and the focus are both determined. On the other hand, the point $-{\bf z}$ is not on the obtained ellipse, since the required ratio for it is $\frac{12-\sqrt{2}}{12}\not =\frac {2}{3}$.

\qed

\begin{theorem}
In a normed plane, an ellipse defined by its leading line and its focus is a convex curve, which is strictly convex if and only if this normed plane is strictly convex.
\end{theorem}

\proof We have to prove that the set
$$
\left \{ {\bf b} \mbox{ : } l\cap \frac{a}{c}\|{\bf b}-{\bf x}\|K=\emptyset \right\}
$$
is a convex domain. It is open by the continuity of the norm function. To prove its convexity, we observe that for any pair of points ${\bf b}_1, {\bf b}_2$ from the considered set the inequality
$$
\frac{a}{c}\|(t{\bf b}_1+(1-t){\bf b}_2)-{\bf x}\|=
\frac{a}{c}\|(t({\bf b}_1-{\bf x})+(1-t)({\bf b}_2-{\bf x})\|\leq $$
$$
\leq t\frac{a}{c}\|({\bf b}_1-{\bf x})\|+(1-t)\frac{a}{c}\|({\bf b}_2-{\bf x})\|
$$
implies that the convex hull of $\frac{a}{c}\|({\bf b}_1-{\bf x})\|K$ and $\frac{a}{c}\|({\bf b}_2-{\bf x})\|K$ contains the ball
$$
\frac{a}{c}\|(t{\bf b}_1+(1-t){\bf b}_2)-{\bf x})\|K.
$$
By the convexity of the half-plane determined by $l$ the convexity property is valid.
Thus the boundary of this domain is a convex curve, as we stated. The second statement is clear from the definition.
\qed

\section{Metric hyperbolas}

\begin{figure}[ht]
    \centering
  \includegraphics[scale=0.5]{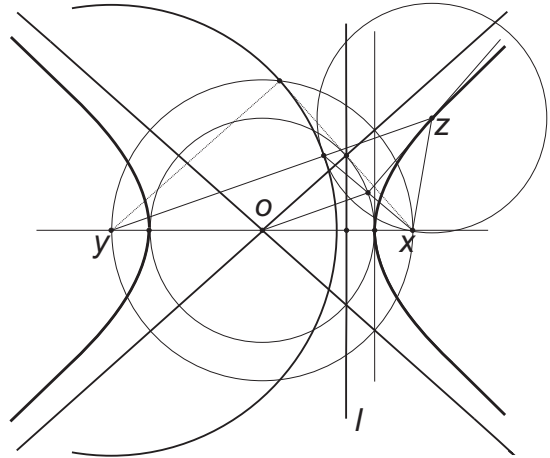}\\
  \caption{The  hyperbola in the Euclidean case}
\end{figure}

     A Euclidean hyperbola satisfies the same metric relations as a Euclidean ellipse, only that now the ratio $\frac{a}{c}$ is smaller than 1. The asymptotes of the hyperbola have directions $\frac{\sqrt{c^2-a^2}}{a}$, and the leading line intersects  the asymptotes in points of the great circle (see Fig. 5). We also have three possible metric definitions. The first one is

\begin{defi}
Given two points {\bf x}, {\bf y} in a normed plane and a distance denoted by $2a > 0$. Then
$$
H({\bf x},{\bf y}, a)=\{ {\bf z} \in X: |\|{\bf z}-{\bf x}\| -\|{\bf y}-{\bf z}\||=2a\}
$$
denotes \emph{the hyperbola defined by its foci} ${\bf x}$ \emph{and} ${\bf y}$. If ${\bf y}=-{\bf x}$, then we use the notation $H({\bf x}, a)$ for it.
\end{defi}

The analogue of Theorem 1 from \cite{wu} is given by our
\begin{theorem}
Let ${\bf x}\in S$ be a point of the unit circle. Then we have:

\noindent (i) $H({\bf x}, 0)$ is the bisector corresponding to the vector {\bf x},

\noindent {(ii)} if there is a neighborhood of {\bf x} on $S$ in which $S$ is strictly convex, then $H(x,2)$ is the union of the two half-lines $[{\bf x},\infty)$ and $[-{\bf x},-\infty)$. If ${\bf x}$ is a point of a piecewise linear part of $S$, then it is the union of two closed cones.
\end{theorem}

\proof
The first statement is obviously true by the definition of the bisector given in the introduction.

The second one follows from the concept and properties of $d$-segments in a Minkowski plane and from our definition of hyperbola; see \cite{martini-swanepoel 1}, \cite{martini-swanepoel 2}, and \S~9 of \cite{B-M-S}.
\qed

From the above theorem it can be seen that a connected part of $H({\bf x},a)$ is, in general, not the boundary of a convex domain, because this property does not hold for a bisector; see \cite{gho1} and \cite{gho2}.

\begin{theorem}
The following two statements are equivalent to each other:

\noindent (i) $K$ is strictly convex.

\noindent (ii) For every ${\bf x} \in S$ and for each value $a \in \mathbb{R}^+$ the set $H({\bf x},a)$ is the union of two simple curves, each of which intersects any line parallel to $[-{\bf x},{\bf x}]$ in precisely two points.

\end{theorem}

\proof From $(i)$ we get $(ii)$. In fact, if $K$ is strictly convex, then every line parallel to $[-{\bf x},{\bf x}]$ contains exactly two points, as said above. This holds by the definition of the hyperbola. These points, dissected by the point of the corresponding bisector, belong to the given line. Thus the sets of the left and right points yield two curves, respectively, which are homeomorphic to the bisector and congruent to each other via reflection in the center of the hyperbola. But the bisector is a simple curve implying the analogous property for the two mentioned curves (see, e.g., \cite{gho1} or \cite{martini-swanepoel 1}).

Conversely, if  $K$ is not strictly convex, then there is a segment in its boundary. Let now ${\bf x}$ be a point of $S$ on the diameter of $K$ parallel to this segment. A hyperbola, corresponding to the foci $\pm{\bf x}$ for every positive $a$ is the union of two closed domains intersected by a line, parallel to ${\bf x}$ and far enough to this diameter, in two segments, since the bisector is also intersected by this line in a seqment. Thus  $(ii)$ does not hold. \qed

\begin{remark}
From this proof we can conclude that the topological properties of hyperbolas do not depend on the parameter $a$ and only on the position of their foci. Thus $(ii)$ is equivalent to

\noindent $(iii)$  \emph{For every ${\bf x} \in S$ there is a value $a \in \mathbb{R}^+\cup \{0\}$ such that the set $H({\bf x},a)$ is the union of two simple curves, intersected by any line parallel to $[-{\bf x},{\bf x}]$ in precisely two points}.
\end{remark}

Analogously to the case of ellipses, we have also two further definitions for hyperbolas. These are given in the following.

\begin{defi}[based on leading circle and focus]
Let $L:=(2a)\cdot K$ be a homothetic copy of the unit disk $K$, and ${\bf x}\in X$ be an arbitrary point exterior to $L$. The locus of points ${\bf z}\in X$ for which there is a positive $\varepsilon$ such that ${\bf z}+\varepsilon K$ touches $L$ and contains ${\bf x }$ on its boundary will be called \emph{the hyperbola defined by its leading circle and its focus} {\bf x}.
\end{defi}

\begin{defi}[based on leading line and focus]
Let $l$ be a straight line, ${\bf x}$ be a point, and $\gamma =\frac{a}{c}$ a ratio less than 1. The locus of points ${\bf z}\in X$, for which there is a positive ${\varepsilon}$ such that the boundary of the disk ${\bf z}+\varepsilon K$ contains {\bf x} and the disk ${\bf z}+\gamma (\varepsilon K)$ touches the line $l$, will be called \emph{the hyperbola defined by its leading line and its focus} {\bf x}.
\end{defi}

Clearly, the three definitions are analogous to those for ellipses.

\begin{proposition}
In normed planes, a hyperbola defined by its foci is always a hyperbola defined by its leading circle and a focus. The converse statement is also true. In general, the third definition yields a different class of curves.
\end{proposition}

\proof  First we consider a hyperbola defined by its leading circle. Then we have two homothetic copies of $K$, $D_{\bf z}={\bf z}+\varepsilon K$ and $L=(2a)K$, in touching position. Then the line joining their centers contains a point from the intersection of their boundaries. Call this point ${\bf p}$. Then
$$
\varepsilon=\|{\bf p}-{\bf z}\|=\|{\bf z}-{\bf x}\|
$$
and
$$
2a=|\|{\bf z}\|-\|{\bf z}-{\bf x}\||,
$$
implying that ${\bf z}\in H({\bf 0},{\bf x}, a)$.

On the other hand, with the same notation we have for a point ${\bf z}\in H({\bf 0},{\bf x}, a)$, by the given definition, that $\|{\bf z}\|>2a$. Consider the disk $D_{\bf z}={\bf z}+(-2a+\|{\bf z}\|)K$. This disk is touching $L$. Since
$$
2a=\|{\bf z}\|-\|{\bf z}-{\bf x}\|
\mbox{ and so }
\|{\bf z}-{\bf x}\|=(-2a+\|{\bf z}\|),
$$
${\bf z}$ is on the boundary of $D_{\bf z}$.

It is clear that the hyperbola defined by its foci is always a centrally symmetric set. This is, in general, not true for a hyperbola defined by its leading line and a focus (see, e.g., the example in Fig. 3). This confirms the final statement in the proposition.
\qed

\begin{figure}[htb]
    \centering
  \includegraphics[scale=0.5]{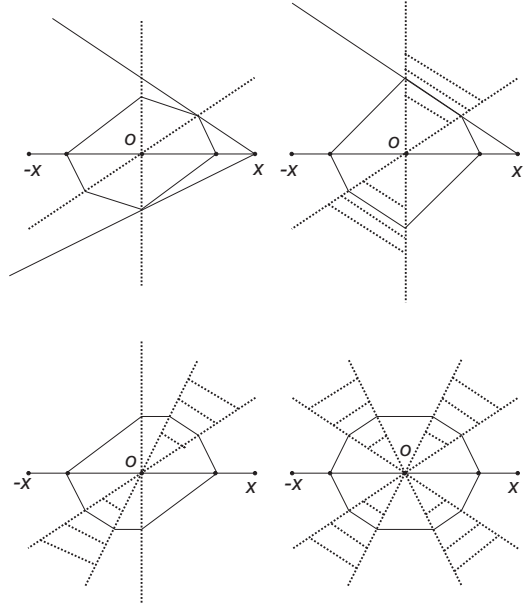}\\
  \caption{The asymptotes}
\end{figure}

\begin{theorem}
The hyperbola defined by its leading line is the union of two simple curves. If the normed plane is strictly convex, then these curves cannot contain segments.
\end{theorem}

\proof We prove that the set $\mathcal{K}$ of exterior points of the hyperbola defined by
$$
\mathcal{K}:=\left\{ {\bf q} \mbox{ : } \frac{a}{c}\|{\bf q}-{\bf x}\|>\rho ({\bf q},l)\right\}
$$
(where $\rho(\cdot,\cdot)$ is the distance function of its arguments) is convex with respect to the direction orthogonal to $l$, where Birkhoff orthogonality is meant. (The concept of \emph{directional convexity} can be found in \cite{gho1}, in connection with the analogous property of the bisector; and  ${\bf x} \in X$ is said to be  \emph{Birkhoff orthogonal} to ${\bf y} \in X$ if $\| {\bf x} + \alpha {\bf y} \| \ge {\bf x}$ holds for any real $\alpha$.) Let {\bf l} be a vector orthogonal to $l$. We will prove that if ${\bf q}$ and ${\bf q}+t{\bf l}$ are in $\mathcal{K}$ for a value $t>0$, then this also holds for every $t'$ with $0\leq t'\leq t$. Now the three points ${\bf q}$, ${\bf q}+t'{\bf l}$, ${\bf q}+t{\bf l}$ are collinear, and their affine hull intersects $l$ in a point {\bf r}. We have distinct cases depending on the position of {\bf r} in this line. For example, assuming the order ${\bf q}, {\bf q}+t{\bf x}, {\bf r}$, we have for all $0\leq t'\leq t$
$$
\rho ({\bf q},l)=\|{\bf r}-{\bf q}\|,
$$
$$
\rho ({\bf q}+t'{\bf l},l)=\|{\bf r}-({\bf q}+t'{\bf l})\|,
$$
and thus
$$
\rho ({\bf q}+t'{\bf l},l)=\rho ({\bf q},l)-t'\|{\bf l}\|.
$$
From this we get that
$$
\frac{a}{c}\|({\bf q}+t'{\bf l})-{\bf x}\|\geq \frac{a}{c}\|{\bf q}-{\bf x}\|-\frac{a}{c}t'\|{\bf l}\|> \rho ({\bf q},l)-t'\|{\bf l}\|=\rho ({\bf q}+t'{\bf l},l)\,,
$$
as we stated. The other cases can be investigated and proved analogously. This means that the boundary of $\mathcal{K}$ is the union of two curves homeomorphic to the line $l$; so they are simple curves intersected by every line parallel to the vector ${\bf l}$ in one point.

Assume that the segment  $[{\bf p}, {\bf q}]$ lies in the hyperbola defined by $l$ and ${\bf x}$. Then we have for all $0\leq t\leq 1$
$$
\frac{a}{c}\|(t{\bf p}+(1-t){\bf q})-{\bf x}\|=\rho (t{\bf p}+(1-t){\bf q},l)=\|(t{\bf p}+(1-t){\bf q})-(t{\bf r}_0+(1-t){\bf r}_1)\|,
$$
where ${\bf r}_t$ is the touching point of the disk having center $t{\bf p}+(1-t){\bf q}$ with its tangent $l$.
Since the disks touching $l$ are homothetic copies of each other, the vectors  ${\bf p}-{\bf r}_0$ and ${\bf q}-{\bf r}_1$ are linearly dependent, and this implies that
$$
\|(t{\bf p}+(1-t){\bf q})-(t{\bf r}_0+(1-t){\bf r}_1)\|=\|t({\bf p}-{\bf r}_0)+(1-t)({\bf q}-{\bf r}_1)\|.
$$
On the other hand,
$$
t\|({\bf p}-{\bf r}_0\|+(1-t)\|{\bf q}-{\bf r}_1\|=t\rho ({\bf p},l)+(1-t)\rho ({\bf q},l)=t\frac{a}{c}\|{\bf p}-{\bf x}\|+(1-t)\frac{a}{c}\|{\bf q}-{\bf x}\|,
$$
and so we also have
$$
\frac{a}{c}\|(t{\bf p}+(1-t){\bf q})-{\bf x}\|=t\frac{a}{c}\|{\bf p}-{\bf x}\|+(1-t)\frac{a}{c}\|{\bf q}-{\bf x}\|
$$
for all $t$. This means that with $t=\frac{1}{2}$ we have
$$
\|({\bf p}+{\bf q})-2{\bf x}\|=\|{\bf p}-{\bf x}\|+\|{\bf q}-{\bf x}\|,
$$
and so, for the points ${\bf \overline{x}} :={\bf p}-{\bf x}, {\bf \overline{y}} :={\bf 0}, {\bf \overline{z}}:=-({\bf q}-{\bf x})$, we get
$$
\|{\bf \overline{x}} -{\bf \overline{z}} \|=\| {\bf \overline{x}} - {\bf \overline{y}}\|+ \| {\bf \overline{y}} -{\bf \overline{z}} \|.
$$
By Proposition 1 in \cite{martini-swanepoel 1} this implies that the segment
$$
\left[\frac{{\bf p}-{\bf x}}{\|{\bf p}-{\bf x}\|}, \frac{{\bf q}-{\bf x}}{\|{\bf q}-{\bf x}\|}\right]
$$
lies in  $S$. Thus the unit ball has to be strictly convex, as we stated.
\qed

The leading circle is a homothetic copy of the main circle with respect to the center {\bf x}. So, if the common tangent of these two circles has unique touching points with these circles, in each case, then there is a natural definition of \emph{asymptotes} joining the center $O$ and the touching points of the tangent lines lying on the main circle. Since asymptotes separate those elements of the pencil of $O$ which intersect the conics from those which are not described by the other case (when the tangent line touches the main circle in a segment), we have the possibility that, as general asymptotes, also conic domains occur.
In Fig. 6 we can see the four possibilities.

\section{Metric parabolas}

\begin{figure}[htb]
  \centering
  \includegraphics[scale=0.7]{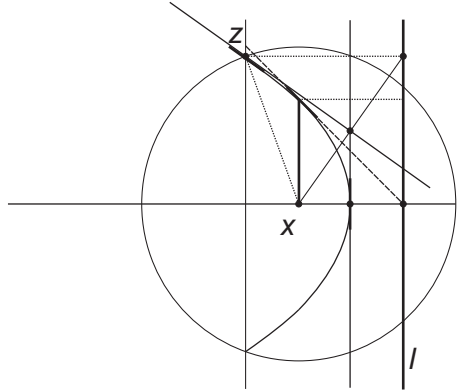}\\
  \caption{The Euclidean case}
\end{figure}

For the case of parabolas, the first two definitions have no analogue, and so we have only the third case.

\begin{defi}
In a normed plane, let $l$ be a straight line, and ${\bf x}$ be a point. The locus of the points ${\bf z}\in S$ for which there is a positive ${\varepsilon}$ such that the boundary of the disk ${\bf z}+\varepsilon K$ contains {\bf x} and touches the line $l$, will be called \emph{the parabola defined by its leading line and its focus} {\bf x}.
\end{defi}

\begin{theorem}
 In a normed plane, the metric parabola is a simple curve which does not contain segments if and only if the normed plane under consideration is strictly convex.
\end{theorem}

\proof
We first prove that if the plane is strictly convex, then any metric parabola is a simple, strictly convex curve, since it is the boundary of a strictly convex domain.

For simplicity, denote by $\rho ({\bf x}, l)$ the radius of that disk whose center is ${\bf x}$ and which touches the line $l$. We first show that the parabola is the common boundary of the two open domains defined by the inequalities
$$
I_{{\bf x},l}:= \|{\bf z}-{\bf x}\|<\rho ({\bf z}, l),
$$
$$
O_{{\bf x},l}:= \|{\bf z}-{\bf x}\|>\rho ({\bf z}, l),
$$
respectively. Consider a point {\bf z} of the parabola.  Then $\|{\bf z}-{\bf x}\|=\rho ({\bf x}, l)$, and by strict convexity there is a unique common point ${\bf l}_{\bf z}$ of $l$ and the circle with center ${\bf z}$ and radius  $\rho ({\bf z}, l)$. The points of the segment $[{\bf x},{\bf z}]$ are in $I_{{\bf x},l}$, and the points of the segments $[{\bf z},{\bf l}_{\bf z}]$ are in $O_{{\bf x},l}$, respectively. This holds since strict convexity implies that if a circle contains another one with smaller radius, then they have at most one common boundary point (see Fig. 8).

\begin{figure}[htb]
  \centering
  \includegraphics[scale=0.5]{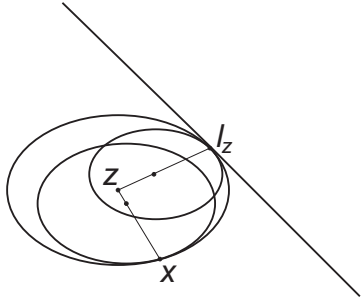}\\
  \caption{The case of strict convexity}
\end{figure}

This means that the point ${\bf z}$ is a common boundary point of the two domains.

Let ${\bf z}_1$ and ${\bf z}_2$ be two points of $I_{{\bf x},l}$. Then the disks with center ${\bf z}_i$ and radius $\rho ({\bf z}_i, l)$ will touch $l$ at the points ${\bf l}_{{\bf z}_i}$, respectively. Observe that if ${\bf z}=t{\bf z}_1+(1-t){\bf z}_2$, where $0\leq t\leq 1$, then
$$
\|{\bf z}-{\bf x}\|\leq t \|{\bf z}_1-{\bf x}\|+(1-t)\|{\bf z}_2-{\bf x}\|<t\rho ({\bf z}_1, l)+(1-t)\rho ({\bf z}_2, l)=t\|{\bf z}_1- {\bf l}_{{\bf z}_1}\|+(1-t)\|{\bf z}_2- {\bf l}_{{\bf z}_2}\|.
$$
But the vectors
${\bf z}_1- {\bf l}_{{\bf z}_1}$ and ${\bf z}_2- {\bf l}_{{\bf z}_2}$ are parallel to each other, and so
$$
t\|{\bf z}_1- {\bf l}_{{\bf z}_1}\|+(1-t)\|{\bf z}_2- {\bf l}_{{\bf z}_2}\|=\|t({\bf z}_1- {\bf l}_{{\bf z}_1})+(1-t)({\bf z}_2- {\bf l}_{{\bf z}_2})\|=\|{\bf z}-{\bf l}_{\bf z}\|,
$$
implying the inequality
$$
\|{\bf z}-{\bf x}\|<\rho ({\bf z}, l)
$$
and the convexity of $I_{{\bf x},l}$. Strict convexity follows from the fact that the inequality in our computation cannot be an equality if the unit disk is strictly convex.

If now the unit disk contains a segment on its boundary, then there is no uniquely determined touching point ${\bf l}_{{\bf z}}$. But we can ``interprete'' all touching segments as ${\bf l}_{{\bf z}}$, and the proof of the first part remains valid. This shows that the parabola is the common boundary of the two domains, $I_{{\bf x},l}$ and $O_{{\bf x},l}$.
The second part of our proof  can also be applied, with the observation that in every situation we can choose points ${\bf y}_i\in {\bf l}_{{\bf z}_i}$, respectively, such that the segments $[{\bf z}_i,{\bf y}_i]$, $i=1,2$, are parallel to each other. Now we can use the first part of the calculation, replacing $<$ by $\leq $, and it is easy to see that by the point
$$
{\bf y}:=t({\bf z}_1- {\bf y}_1)+(1-t)({\bf z}_2- {\bf y}_2)
$$
we can represent the length $\rho({\bf z},l)$, again confirming convexity. Now our proof is complete if we observe that in case of  $l_{\infty}$ there is a metric parabola containing segments, as we can also see in Fig. 3.
\qed

\section{Some further remarks on conics}

Of course, the projective geometry of a normed space and its embedded Euclidean geometry is the same. But there is a possibility to take into consideration the metric, too, because in a ``nice'' normed space there is a theory of selfadjoint linear transformations.

The following way is a possibility to define quadrics in the projective augmentation of any smooth, strictly convex space. We describe this method in the two-dimensional case, where the quadric is clearly a conic.

Every normed plane can be represented as a semi-inner product space (s.i.p.; see \cite{lumer} and \cite{giles}). If the unit disk  is strictly convex, this representation is unique. As proved in \cite{giles}, the orthogonality with respect to the s.i.p. is equivalent to the orthogonality concept of Birkhoff (see, e.g., \cite{alonso1} and \cite{alonso2}).  Koehler proved in \cite{koehler} that if the generalized Riesz-Fischer
representation theorem is valid in a normed space, then every bounded linear operator $A$ has a generalized adjoint $A^T$ defined by the equality
$$
<A(x),y> = <x,A^T(y)> \mbox{ for all } x,y \in V.
$$
It can be proved that if in all strictly convex and smooth spaces the above assumption holds, then in such a space there is a generalized adjoint.
We remark that $A^T$ is in general not a linear transformation. We say that the linear mapping is \emph{self-adjoint} if $A=A^T$.
If $A$ is self-adjoint, then any element of its class in the Projective General Linear Group of $V$ is self-adjoint, too. So we can call such a family of
operators \emph{class of self-adjoint linear operators of the projective space} $P(V)$. Now the concept of conics can be introduced as follows.

\begin{defi}
Let $P(V)$ be a real projective space with the two-dimensional semi-inner product space $(V,<\cdot,\cdot>)$. A (non-degenerate) projective conic is the zero set of a (non-degenerate) form $\Phi(x, y) = <A(x),y>$, with an invertible self-adjoint operator $A$ of $P(V)$.
\end{defi}

We remark that the form $\Phi(x,y)$ is linear in its first argument, homogeneous in its second one, but is neither symmetric, bilinear nor positive. It is symmetric and bilinear if the semi-inner product is symmetric; bilinear if the semi-inner product is additive in its second argument; and positive if $A$ is a square operator (meaning that it is the square of another self-adjoint operator, denoted by $\sqrt{A}$).

The group of self-adjoint operators is basically determined by the unit disks, and it determines the projective conics. Thus, in this
setting the metric of the plane is also used for smooth, strictly convex
normed planes.

We finish with two problems:
\begin{enumerate}
\item Characterize the self-adjoint operators for smooth, strictly convex normed planes.
\item Describe relations between metric conics and general ones.
\end{enumerate}

The first question was also investigated in \cite{langi} in the case when the plane also has a Lipshitz-type property. The second question is the theme of a forthcoming paper of the authors.

\end{document}